\documentclass[12pt,a4paper]{article}

\usepackage{latexsym,amssymb,amsthm}
\usepackage{url,graphics}
\usepackage[dvips]{graphicx,epsfig}

\title{\bf On existence of noncritical vertices in digraphs}
\author{
     \begin{tabular}{c}
     G.\,V.\,Nenashev
     \\[6pt]
     \\[-2pt]
     {\small E-mail: \texttt{glebnen@mail.ru}       }
    \end{tabular}
    }
\date{}

\begin{document}
\maketitle
\righthyphenmin=2
\renewcommand*{\proofname}{\bf Proof}
\newtheorem{thm}{Theorem}
\newtheorem{lem}{Lemma}
\newtheorem{cor}{Corollary}
\theoremstyle{definition}
\newtheorem{defin}{Definition}
\theoremstyle{remark}
\newtheorem{rem}{\bf Remark}

\def\I{{\rm Int}}
\def\R{{\rm Bound}}
\def\q#1.{{\bf #1.}}
\renewcommand\geq{\geqslant}
\renewcommand\leq{\leqslant}

\abstract{ Let~$D$ be a strongly connected digraphs on~$n\ge 4$ vertices. A vertex~$v$ of~$D$ is noncritical, if the digraph~$D-v$ is strongly connected. We  prove, that if sum of the degrees of any two adjacent vertices of~$D$ is at least $n+1$, then there exists a noncritical vertex in~$D$, and if sum of the degrees of any two adjacent vertices of~$D$ is at least $n+2$, then there exist two noncritical vertices in~$D$.  A series of examples confirm that these bounds are tight.}

\section{\bf Introduction}
In this paper we consider a digraph~$D$ without loops and multiple arcs, i.e. any two vertices~$x$ and~$y$ are connected by at most two arcs
(at most one arc of each direction).

For a subgraph~$T$ of the digraph~$D$ we denote by~$V(T)$ the set of vertices of~$T$ and by~$\overline{V(T)}$ the set of vertices of~$D$ which do not belong to the subgraph~$T$.
We denote by~$T-x$  the subgraph  obtained from~$T$ by deleting the vertex~$x$ and all  arcs incident to~$x$.

\begin{defin}
A digraph is called {\it strongly connected} if for any two its vertices~$x$, $y$ there is a path from~$x$ to~$y$.
\end{defin}

\begin{defin}
A vertex~$v$ of a strongly connected digraph~$D$ is called  {\it noncritical}, if the digraph~$D-v$ is strongly connected.
\end{defin}

\begin{defin}
{\it The degree} of a vertex~$x$ in a digraph~$D$ (notation: $deg(x)$) is the number of vertices adjacent to~$x$. (In the case where vertices~$x$ and~$y$ are connected by several arcs we count the vertex~$y$ once.)

\end{defin}

S.\,V.\,Savchenko~\cite{S1} in 2006 has proved that a strongly connected digraph~$D$ with~$n$ vertices and vertex degrees at least~$\frac{3n}{4}$ has two noncritical vertices. There are no further  results  on this subject. 

It follows from our main result, that if minimal vertex degree of a digraph  is at least $\frac{n+1}{2}$ then this digraph has a noncritical vertex, and if minimal vertex degree of a digraph  is at least $\frac{n+1}{2}$ then this digraph has two noncritical vertices.

As well as the author of~\cite{S1}, we need the following lemma, formulated in~\cite{Ma}.

\begin{lem}
Let $D$ be a strongly connected digraph and~$S$ be a proper strongly connected subgraph of~$D$. Then~$S$ is a maximal proper strongly connected subgraph of~$D$ if and only if the three following conditions hold:

{\bf 1)} There exists a vertex~$\omega_{in} \in \overline{V(S)}$, such that any arc from~$V(S)$ to~$\overline{V(S)}$ ends in the vertex~$\omega_{in}$;
\smallskip

{\bf 2)} There exists a vertex~$\omega_{out} \in \overline{V(S)}$, such that any arc from~$\overline{V(S)}$ to~$V(S)$ begins at the vertex~$\omega_{out}$;
\smallskip

{\bf 3)} There exists a unique simple path from~$\omega_{in}$ to~$\omega_{out}$ in~$D$  and this path contains all  vertices of the set~$\overline{V(S)}$ and only them (see figure~$1$).
\end{lem}

\section{\bf Search for noncritical vertices}

\begin{thm}
\label{t1} Let~$D$ be a strongly connected digraph  with $n\geq 4$ vertices, such that for any two adjacent vertices~$x$ and~$y$ of this digraph  the inequality~$deg(x)+deg(y)\geq n+1$ holds. Then~$D$ has a noncritical vertex.

\begin{proof}
We consider two cases~$A$ and $B$: in case~$A$ the minimal vertex degree of~$D$ is at least~$3$, and in case~$B$ there is a vertex of degree less than~$3$ in~$D$.

{\bf A.} {\it The minimal vertex degree of~$D$ is at least~$3$.}

Consider a maximal proper strongly connected subgraph~$S$ of the digraph~$D$.

${\bf 1}^\circ $ {\it If $\overline{V(S)}$ consists of one vertex}, then this vertex is noncritical and the theorem is proved.
\smallskip

Hence in what follows the set~$\overline{V(S)}$ consists of at least two vertices, then the vertices~$\omega_{in}$ and~$\omega_{out}$ are different.
\smallskip

${\bf 2}^\circ$ {\it Let $\overline{V(S)}$ contains at least four vertices}.

Consider a path from~$\omega_{in}$ to~$\omega_{out}$, which contains all vertices of the set~$\overline{V(S)}$ and only them.
In this case we can choose  two successive vertices~$a_1$ and~$a_2$ in this path, which are different from~$\omega_{in}$ and~$\omega_{out}$ (see figure~1).

{\begin{figure}[htb!]
\centering
\includegraphics[scale=0.4]{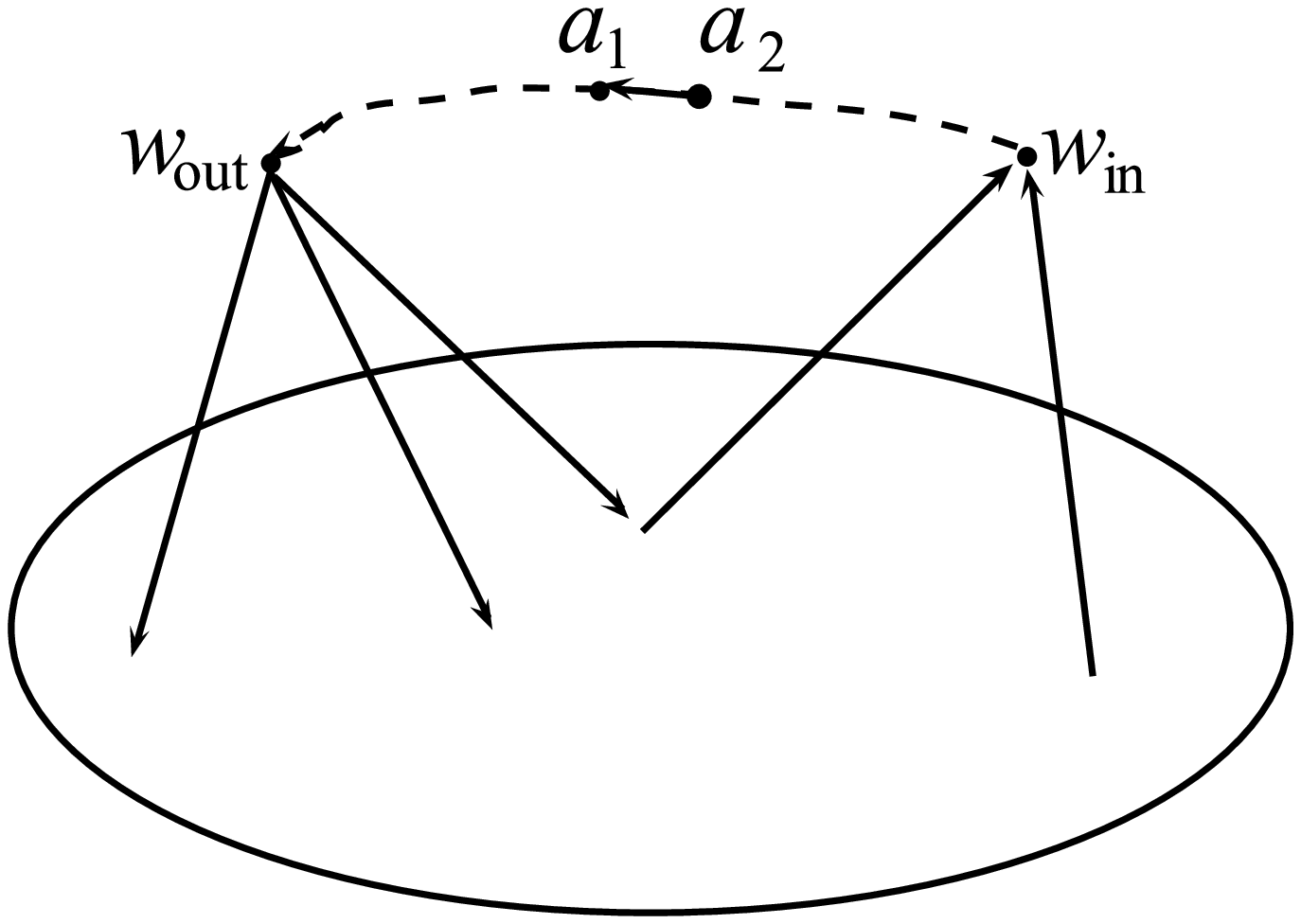}
\caption{ }
\end{figure}
}

Then there are no arcs between~$\{a_1,a_2\}$ and~$S$, hence the degree of both  vertices~$a_1$ and~$a_2$  is at most~$|\overline{V(S)}|-1$.
Consider two adjacent vertices~$b_1, b_2 \in V(S)$.  (These vertices exist since otherwise~$S$ consists of one vertex and this vertex has degree at most~2.) The degree of each of vertices~$b_1$, $b_2$ is not more than~$|V(S)|-1+2=|V(S)|+1$ (the vertex~$b_i$ can be adjacent to 
the vertices of~$S$ different from~$b_i$ and to~$\omega_{in}$, $\omega_{out}$).

We have chosen two pairs of adjacent vertices, hence, we obtain the following inequality:
$$2(n+1) \leq deg(a_1)+deg(a_2)+deg(b_1)+deg(b_2) \leq $$
$$2(|\overline{V(S)}|-1) + 2(|V(S)|+1) = 2|V(D)| = 2n, $$
that is impossible.

\smallskip

${\bf 3}^\circ$ {\it Let $\overline{V(S)}$ contains three vertices}. 

In this case a vertex~$a\in \overline{V(S)}$, which is different from~$\omega_{in}$ and~$\omega_{out}$,  is adjacent only to $\omega_{in}$  and~$\omega_{out}$, i.e.~$deg(a)=2$. We obtain a contradiction.
\smallskip

${\bf 4}^\circ${\it Let $\overline{V(S)}$ contains two vertices},
i.e.~$\overline{V(S)}=\{\omega_{in},\omega_{out}\}$.

{\begin{figure}[hbt!]
\centering
\includegraphics[scale=0.4]{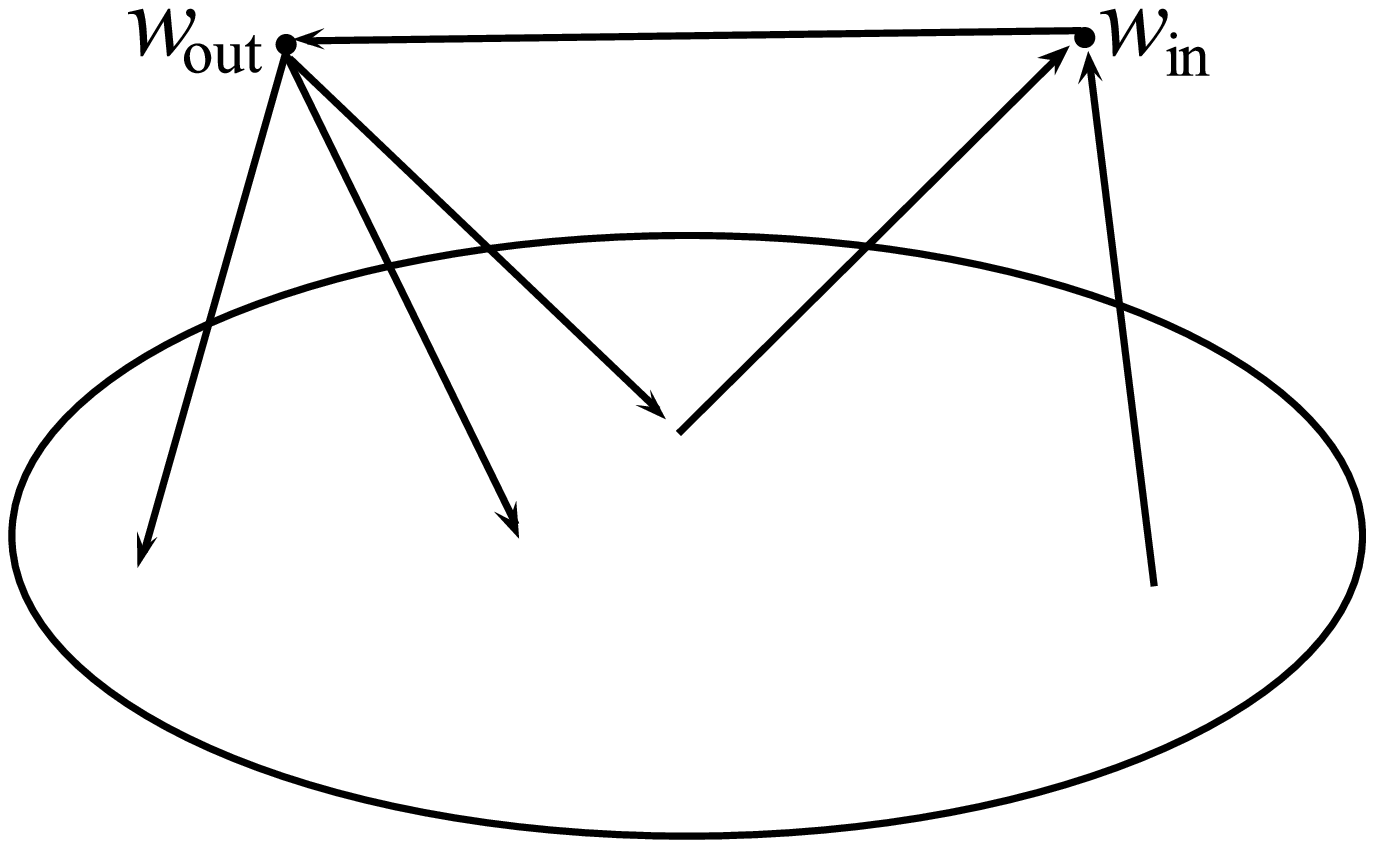}
\caption{ }
\end{figure}
}

Let us construct an {\it incoming tree} $T_{in}$ of the vertex~$\omega_{in}$ in the subgraph~$D-\omega_{out}$. The vertex~$\omega_{in}$ is {\it the root} of this tree. Level~1 consists of the vertices which have  outgoing  arcs to~$\omega_{in}$, and so on:
level~$k$ consists of the vertices which do not belong to previous levels and have outgoing  arcs to vertices of level~$(k-1)$. 
It follows from the strong connectivity of~$D$, that for any vertex~$x\in S$  there is a path from~$x$ to~$\omega_{in}$ in~$D$. Hence, any vertex of the set~$S$ belongs to some level. For each vertex, with the exception  of~$\omega_{in}$, we draw exactly one outgoing arc to a vertex of previous level (see figure~3).

{\begin{figure}[htb!]
\centering
\includegraphics[scale=0.4]{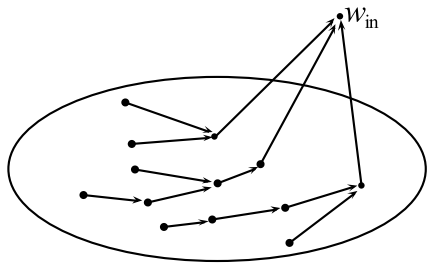}
\caption{Incoming tree.}
\end{figure}
}

Now consider~$T_{in}$ as an undirected tree. Let us call any connected component of~${T_{in} - \omega_{in}}$ a {\it branch}. Note, that every branch contains a leaf of this tree, and the number of branches is equal to~$deg(\omega_{in})-1$, since~$T_{in}$ contains all arcs incident  to~$\omega_{in}$,  with the exception of  arcs between~$\omega_{in}$ and~$\omega_{out}$.
Denote by~$A_{in}$ the set of all leaves of~$T_{in}$. We have proved that~$|A_{in}|\geq deg(\omega_{in})-1$.
Clearly, for any vertex~$x\in A_{in}$, there is a path from any other vertex of the set~$S$ to~$\omega_{in}$ in~$T_{in}-x$.

Similarly, we  construct an {\it outgoing tree} of the vertex~$\omega_{out}$ in the subgraph~$D-\omega_{in}$ and the set~$A_{out}$.
For any vertex~$x\in A_{out}$, there is a path from~$\omega_{out}$ to any other vertex of the set~$S$  in~$T_{out}-x$.
Similarly,~$|A_{out}|\ge deg(\omega_{out})-1$. The vertices~$\omega_{in}$ and~$\omega_{out}$ are adjacent, whence it follows:
$$ |A_{in}|+|A_{out}|\geq deg(\omega_{in})-1+deg(\omega_{out})-1\geq n+1-2= n-1.$$
Since $|A_{in}\cup A_{out}| \le |S| =n-2$, there exists a vertex~$x\in A_{in}\cap A_{out}$.

For any vertex~$y\in S\setminus x$ there is a path from~$y$ to~$\omega_{in}$ and a path from~$\omega_{out}$ to~$y$ in the graph~$D-x$. Since there is an arc from~$\omega_{in}$ to~$\omega_{out}$, the graph~$D-x$  is strongly connected and the vertex~$x$ is noncritical.

\smallskip
{\bf B.} {\it There is a vertex of degree less than~$3$ in~$D$.} 

Let~$q$ is a vertex of degree less than~$3$. Clearly, there exists a vertex~$p_1$ adjacent to~$q$. 
We know, that~$$deg(q)+n-1 \geq deg(q)+deg(p_1)\geq n+1,$$
 hence we obtain~$deg(q)=2$ and~$deg(p_1)=n-1$. Thus there exists another vertex~$p_2$ adjacent to~$q$.
Similarly, $deg(p_2)=n-1$.

Since~$deg(p_1)=deg(p_2)=n-1$, there exists an arc between~$p_1$ and~$p_2$ (maybe two arcs of different directions). Without loss of generality assume that there is an arc~$p_1p_2$ in~$D$.  We assume that~$q$ is not a noncritical vertex (otherwise the theorem is proved). Then the graph~$D-q$ is not strongly connected, hence there are two vertices~$x$ and~$y$, such that any path from~$x$ to~$y$ in~$D$ contains the vertex~$q$. Consider the shortest path~$P$ from~$x$ to~$y$. Clearly,~$P$ must pass both vertices~$p_1$ and~$p_2$. Since~$P$ is the shortest path, it passes the arcs~$p_2q$ and~$qp_1$ (otherwise the shortest path must  pass the arc~$p_1p_2$ and avoid~$q$).

Hence there is an oriented cycle~$qp_1p_2$ in the digraph~$D$. Since this cycle does not contain all vertices of~$D$, there is a maximal strongly connected proper subgraph~$S$, which contains~$q$, $p_1$  and~$p_2$. If ${|\overline{V(S)}|>}2$, then there exists a vertex~$a\in \overline{S}$, different from~$\omega_{in}, \omega_{out}$. Clearly, there is no arc between~$p_1$ and~$a$, i.e. $deg(p_1)<n-1$, we obtain a contradiction.
The remaining cases~$|\overline{V(S)}|=1$ and~$|\overline{V(S)}|=2$  are similar to the subcases~$1^\circ$ and~$4^\circ$ of  case~$A$. 
(It does not matter whether there are vertices of degrees~1 and~2 or not in the proofs of these subcases).
\end{proof}

\end{thm}

\begin{cor}
\label{cor11}
Let~$D$ be a strongly connected digraph with~$n\geq 4$ vertices and vertex degrees at least~$\frac{n+1}{2}$. Then~$D$ has a noncritical vertex.
\end{cor}

\begin{rem}
1) The bounds~$n+1$ in theorem~\ref{t1} and~$\frac{n+1}{2}$ in corollary~\ref{cor11} are tight. Let us construct  for an even~$n$ a graph~$D$ (see figure~4), such that:

--- $V(D)=\{a_1,\ldots a_{n/2}$,$b_1,\ldots b_{n/2}\}$;

--- $E(D)$  consists of arcs~$a_ib_i$ and arcs~$b_ia_j$, where~$i\neq j$.

Clearly, all vertex degrees in this graph are equal to~$n\over2$ and there are no noncritical vertices.
{\begin{figure}[htb!]
\centering
\includegraphics[scale=0.4]{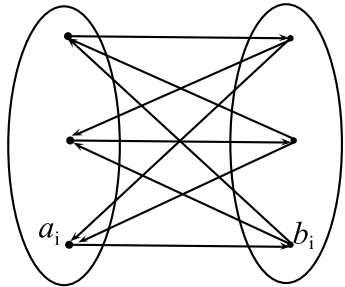}
\caption{A graph without noncritical vertices. }
\end{figure}
}\\

2)  The bound~$n+1$ in theorem~\ref{t1}  is also tight for odd~$n$. 
Let us construct a digraph~$D$ which is suitable for all~$n\geq 6$ (see figure~5):

--- $V(D)=\{a_1,\ldots, a_{n-4}$,$x_1,\ldots, x_4\}$.

--- $E(D)$  consists of arcs~$a_ia_{i+1}$, $a_ja_i$ ($j>i+1$), $a_ix_3$, $x_2a_i$, $a_{n-4}x_1$, $x_1x_2$, $x_2x_3$, 
 $x_3x_4$,  $x_4a_1$ (where~$i,j\in\{1,\dots, n-4\}$).

{\begin{figure}[htb!]
\centering
\includegraphics[scale=0.75]{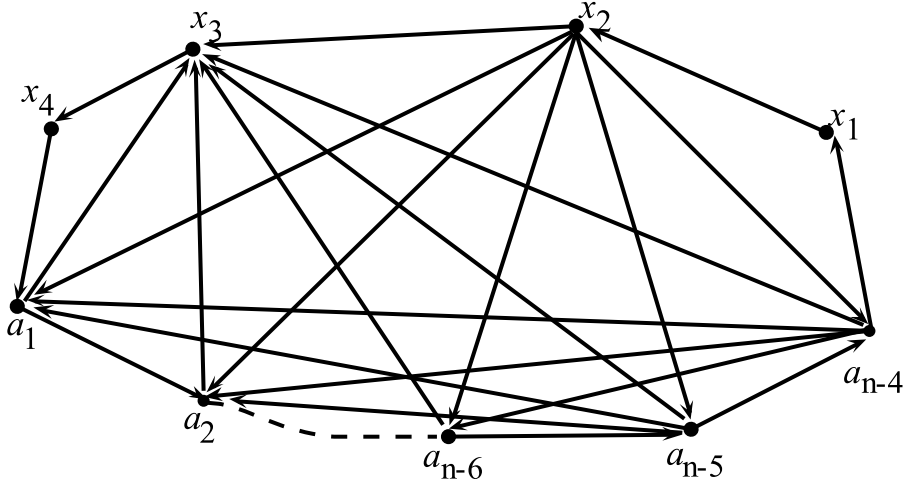}
\caption{A graph without noncritical vertices.}
\end{figure}
}

Clearly, the digraph~$D$ is strongly connected. Let us verify the condition on sum of degrees of pairs of adjacent vertices. 
The vertices~$x_1$ and~$x_4$ have degree~2, the vertices $a_1,a_{n-4}$, $x_2,x_3$ have degree~$n-2$, all other vertices have degree~$n-3$.
Each vertex of  degree~2 is adjacent only to vertices of degree~$n-2$, i.e.  the sum of degrees of these pairs is~$n$. 
Any other vertex has degree at least~$n-3$, hence, any other pair of adjacent vertices has sum of degrees at least~$2(n-3)=2n-6\geq n$.

Let us assure, that~$D$ has no noncritical vertices. In the graph~$D-a_i$ (where~$i>1$) there is no path from~$\{x_3,x_4,a_1,\ldots, a_{i-1}\}$ to~$x_1$. In the graph~$D-x_1$ there are no arcs with the end at~$x_2$. In the graph~$D-x_2$ there are no arcs with the beginning at~$x_2$.   In the graph~$D-x_3$ there are no arcs with the end at~$x_4$. In the graph~$D-x_4$ there are no arcs with the beginning at~$x_3$.
In the graph~$D-a_1$ there are no arcs with the beginning at~$x_4$. Hence, there are no noncritical vertices in~$D$.
\end{rem}

\begin{cor}
\label{cor21}
Let~$D$ be a strongly connected digraph  with $n\geq 4$ vertices, such that for any two adjacent vertices~$x$ and~$y$ of this digraph  the inequality~$deg(x)+deg(y)\geq n+2$ holds. Then~$D$ has two noncritical vertices.

\begin{proof}
We claim, that the degree of any vertex  in the graph~$D$ is at least~3. (Otherwise, if there is a vertex~$x$ of  degree at most two, then  a vertex adjacent to~$x$ must have degree at least~$n$. Clearly, that is impossible). Hence, we have the case~$A$ of theorem~1. Consider two cases.

${\bf 1}^\circ$ {\it For any maximal strongly connected proper subgraph~$S$ we have $|\overline{V(S)}|=1$.}

By theorem~\ref{t1} there exists a noncritical vertex~$x_1$ in~$D$. Consider a maximal strongly connected proper subgraph which contains~$x_1$.  Let it does not contain~$x_2$. Then~$x_1$ and~$x_2$ are two different noncritical vertices.
\smallskip

${\bf 2}^\circ$ {\it There exists a maximal strongly connected subgraph~$S$, such that~$|\overline{V(S)}|\geq2$.}

Then by the reasonings of theorem~\ref{t1}   only the  case where ${|\overline{V(S)}|=2}$ is possible.
Consider the sets~$A_{in}$ and~$A_{out}$ constructed in the proof of theorem~\ref{t1} (see subcase~$4^\circ$ of case~$A$). Now we have the inequality
$$ |A_{in}|+|A_{out}|\geq deg(\omega_{in})-1+deg(\omega_{out})-1\geq n+2-2= n.$$
Since $|A_{in}\cup A_{out}| \le |S| =n-2$, there exist at least two  vertices in~$A_{in}\cap A_{out}$.
Clearly, these vertices are noncritical.
\end{proof}
\end{cor}

\begin{cor}
\label{cor22}
Let~$D$ be a strongly connected digraph with~$n\geq 4$ vertices and vertex degrees at least~$\frac{n+2}{2}$. Then~$D$ has at least two noncritical vertices.
\end{cor}

\begin{rem}

1) The bounds~$n+2$ in corollary~\ref{cor21} and~$\frac{n+2}{2}$ in corollary~\ref{cor22} are tight. Let us construct  for an odd~$n$ a graph~$D$ (see figure~6), such that:

--- $V(D)=\{a_1,\ldots a_{(n-1)/2},b_1,\ldots b_{(n-1)/2},x\}$;

--- $E(D)$  consists of arcs~$a_ib_i$, arcs~$b_ia_j$, where~$i\neq j$,  arcs~$xa_i$ and arcs~$b_ix$.

Clearly, all vertex degrees in this graph are at least~$n+1\over2$ and~$x$ is the only  noncritical vertex.
{\begin{figure}[hbt!]
\centering
\includegraphics[scale=0.4]{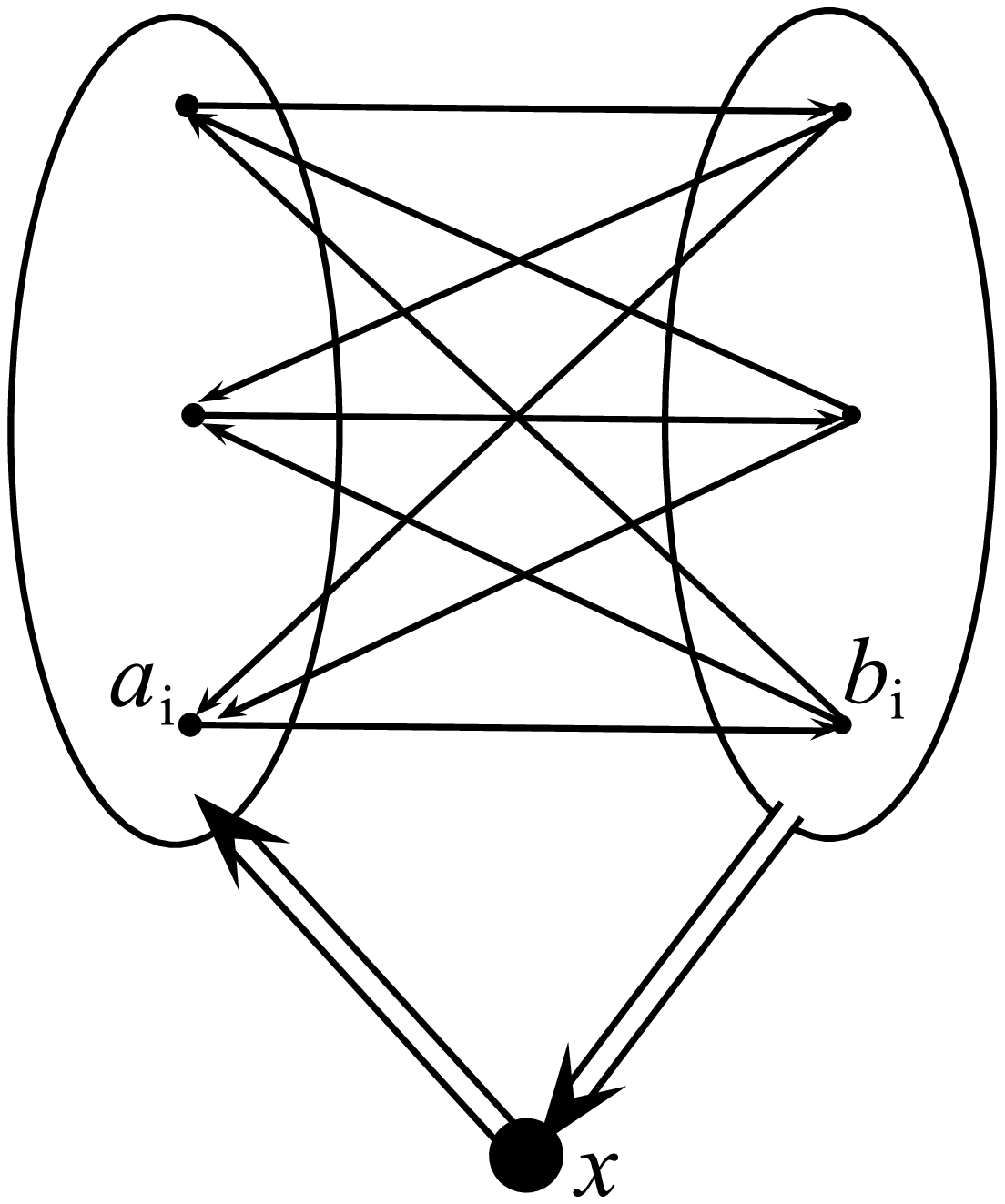}
\caption{A graph with only one noncritical vertex.}
\end{figure}
} \\

 2)  The bound~$n+2$ in corollary~\ref{cor21} is also tight for even~$n$. 
Let us construct a digraph~$D$ which is suitable for all~$n\geq 5$ (see figure~7):

--- $V(D)=\{a_1,\ldots, a_{n-1}$,$x\}$.

--- $E(D)$ consists of arcs~$a_ia_{i+1}$, $a_ja_i$ ($j>i+1$), $xa_1$, $a_{n-1}x$ (where~$i,j\in\{1,\dots, n-1\}$).

{\begin{figure}[htb!]
\centering
\includegraphics[scale=0.75]{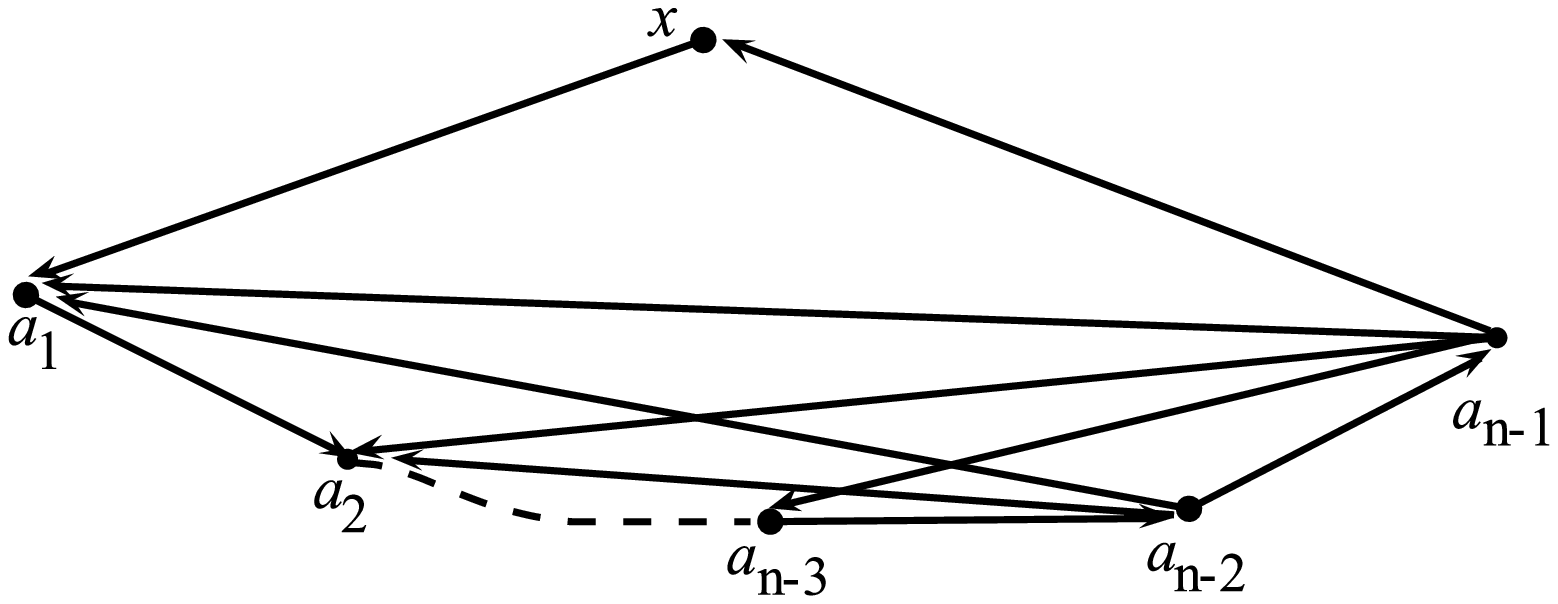}
\caption{A graph with only one noncritical vertex.}
\end{figure}
}

Clearly, the digraph~$D$ is strongly connected. Let us verify the condition on sum of degrees of pairs of adjacent vertices. 
The vertex~$x$ of degree~2 is adjacent in~$D$ only to vertices of degree~$n-1$, i.e.  the sum of degrees of these pairs is~$n+1$. 
Any other vertex has degree at least~$n-2$, hence, any other pair of adjacent vertices has sum of degrees at least~$2(n-2)=2n-4\geq n+1$.

Let us assure, that~$x$ is the only noncritical vertex in~$D$. Clearly, in the graph~$D-a_i$ (where~$i>1$) there is no path from~$\{a_1,\ldots, a_{i-1}\}$ to~$x$.  In the graph~$D-a_1$ there are no arcs with beginning at~$x$.

\end{rem}

\begin{rem}
The proof of theorem~\ref{t1} in the case~$|\overline{V(S)}|>1$ and minimal vertex degree is~$d$ gives us~$2d-n$ noncritical vertices. 
In spite of this, no lower bound on minimal vertex degree can provide three noncritical vertices in a strongly connected digraph~$D$. 

For any~$n$ there exists a strong tournament  (i.e. a strongly connected digraph in which each pair of vertices is connected by exactly one arc) with precisely two noncritical vertices. Let us construct this tournament on~$n$ vertices~$a_1,\dots, a_n$  (see figure~8). This tournament has arcs of type~$a_ia_{i+1}$ (where~$i<n$), and of type~$a_ia_j$ for each pair~$i,j$, where~$i>j+1$.

{\begin{figure}[htb!]
\centering
\includegraphics[scale=0.75]{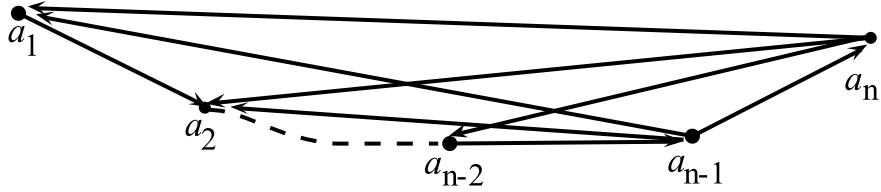}
\caption{Strong tournament with two noncritical vertices.}
\end{figure}
}

\end{rem}

\smallskip Translated by D.\,V.\,Karpov.

\end{document}